\documentclass[envcountsame, envcountsect]{llncs}

\usepackage{amssymb,amsmath}

\numberwithin{equation}{section}

\newcommand{\F}{\mathbb{F}}
\newcommand{\Q}{\mathbb{Q}}
\newcommand{\Z}{\mathbb{Z}}

\begin{document}

\title{Constructing Pairing-Friendly Elliptic Curves with Embedding Degree $10$.}
\author{David Freeman}
\institute{University of California, Berkeley \\
	{\tt dfreeman@math.berkeley.edu}}
\maketitle

\begin{abstract}
    We present a general framework for constructing families of elliptic curves of prime order with prescribed embedding degree.  We demonstrate this method by constructing curves with embedding degree $k = 10$, which solves an open problem posed by Boneh, Lynn, and Shacham \cite{bls}.  We show that our framework incorporates existing constructions for $k = 3$, $4$, $6$, and $12$, and we give evidence that the method is unlikely to produce infinite families of curves with embedding degree $k > 12$.
\end{abstract}

\section{Introduction}
\label{s:intro}

A cryptographic pairing is a bilinear map between two groups in which the discrete logarithm problem is hard.  In recent years, such pairings have been applied to a host of previously unsolved problems in cryptography, the most important of which are one-round three-way key exchange \cite{joux}, identity-based encryption \cite{bf}, and short digital signatures \cite{bls}.

The cryptographic pairings used to construct these systems in practice are based on the Weil and Tate pairings on elliptic curves over finite fields.  These pairings are bilinear maps from an elliptic curve group $E(\F_q)$ to the multiplicative group of some extension field $\F_{q^k}$.  The parameter $k$ is called the {\it embedding degree} of the elliptic curve.  The pairing is considered to be secure if taking discrete logarithms in the groups $E(\F_q)$ and $\F_{q^k}^*$ are both computationally infeasible.  

For optimal performance, the parameters $q$ and $k$ should be chosen so that the two discrete logarithm problems are of approximately equal difficulty when using the best known algorithms, and the order of the group $\#E(\F_q)$ should have a large prime factor $r$.  For example, a pairing is considered secure against today's best attacks when $r \sim 2^{160}$ and $k \sim 6 \mbox{-} 10$, depending on the application.  In order to vary the security level or adapt to future improvements in discrete log technology, we would like to have a supply of elliptic curves at our disposal for arbitrary $q$ and $k$.

Many researchers have examined the problem of constructing elliptic curves with prescribed embedding degree.  Menezes, Okamoto, and Vanstone \cite{mov} showed that a supersingular elliptic curve must have embedding degree $k \leq 6$, and furthermore $k \leq 3$ in characteristic not equal to $2$ or $3$.  Miyaji, Nakabayashi, and Takano \cite{mnt} have given a complete characterization of ordinary elliptic curves of prime order with embedding degree $k = 3$, $4$, or $6$, while Barreto and Naehrig \cite{bn} give a construction for curves of prime order with $k = 12$.  There is a general construction, originally due to Cocks and Pinch \cite{cp}, for curves of arbitrary embedding degree $k$, but in this construction the sizes of the field $\F_q$ and the subgroup of prime order $r$ are related by $q \approx r^2$, which leads to inefficient implementation.  Recent efforts (cf. \cite{bw}, \cite{cdc}) have focused on reducing the ratio $\rho = \log q/\log r$ for arbitrary $k$, but no additional examples have been found with $\rho$ small enough to allow for curves of prime order.

The focus of this paper is the construction of ordinary elliptic curves of prime order with prescribed embedding degree.  In Section \ref{s:method} we present a general framework for constructing such curves and give conditions under which this method will give us infinite families of elliptic curves.  The method is based on the Complex Multiplication method of curve construction \cite{mor} and is implicit in the constructions of several other researchersl.  Our contribution is to gather all of the relevant results in one place and to define terminology that makes it apparent that these various constructions are all instances of the same general method.

Our main contribution appears in Section \ref{s:k=10}, where we show how the method of Section \ref{s:method} can be used to construct curves with embedding degree $k = 10$.  We give examples of such curves over fields of cryptographic size, solving an open problem posed by Boneh, Lynn, and Shacham \cite{bls}. 
 
In Section \ref{s:small} we show how the existing constructions of elliptic curves of prime order with embedding degree $k = 3$, $4$, $6$, or $12$ can be explained via the  framework of Section \ref{s:method}.  In Section \ref{s:higher}, we show that for $k > 6$, our method is not likely to give additional infinite families of elliptic curves with the specified embedding degree.  We note, however, that examples of such families exist for $k = 10$ and $k = 12$, and we ask in Section \ref{s:conclusion} if such examples can be constructed in a systematic fashion.

\subsection*{Acknowledgments}

Research for this paper was conducted during a summer internship at Hewlett-Packard Laboratories, Palo Alto.  I thank Vinay Deolalikar for suggesting this topic and for providing advice and support along the way.  I also thank Gadiel Seroussi for bringing me to HP and for supporting my research.

I thank Paulo Barreto, Steven Galbraith, Ed Schaefer, and Mike Scott for their valuable feedback on earlier versions of this paper.  I am especially indebted to Mike Scott, who used the method presented in Section \ref{s:k=10} to compute examples of elliptic curves of cryptographic size with embedding degree $10$.  Two of these curves now appear in this paper as Examples \ref{ex1} and \ref{ex2}.

\section{A Framework for Constructing Pairing-Friendly Elliptic Curves}
\label{s:method}

In this section we describe a general framework for constructing elliptic curves of a given embedding degree $k$.  This framework is implicit in the constructions of Miyaji, Nakabayashi, and Takano \cite{mnt}; Barreto, Lynn, and Scott \cite{bls1}; Cocks and Pinch \cite{cp} (as explained in \cite{bss2}); and Brezing and Weng \cite{bw}.  After stating the relevant results, we define terminology that will allow us to show that these constructions are all specific cases of the same general method.

To construct our elliptic curves, we parameterize the number of points on the curve and the size of the field of definition by polynomials $n(x)$ and $q(x)$, respectively.  For each $x_0$ that gives prime values for $n(x_0)$ and $q(x_0)$, we can use the Complex Multiplication method to construct an elliptic curve with the desired properties.   The main result of this section is Theorem \ref{t:quadratic}, which gives a criterion for the existence of infinite families of such good parameters. 

We begin by giving a formal definition of embedding degree.

\begin{definition}
Let $E$ be an elliptic curve defined over a finite field $\F_q$, and let $n$ be a prime dividing $\#E(\F_q)$.  The {\it embedding degree of $E$ with respect to $n$} is the smallest integer $k$ such that $n$ divides $q^k-1$.
\end{definition}

Equivalently, $k$ is the smallest integer such that $\F_{q^k}$ contains $\mu_n$, the group of $n$th roots of unity in $\overline{\F}_q$.  We often ignore $n$ when stating the embedding degree, as it is usually clear from the context.

If we fix a target embedding degree $k$, we wish to solve the following problem: find a prime (power) $q$ and an elliptic curve $E$ defined over $\F_q$ such that $n = \#E(\F_q)$ is prime and $E$ has embedding degree $k$.  Furthermore, since we may wish to construct curves over fields of different sizes, we would like to be able to specify (approximately) the number of bits of $q$ in advance.

We follow the strategy of Barreto and Naehrig \cite{bn} in parameterizing the trace of the curves to be constructed. Namely, we choose some polynomial $t(x)$, which will be the trace of Frobenius for our hypothetical curve, and construct polynomials $q(x)$ and $n(x)$ that are possible orders of the prime field and the elliptic curve group, respectively.  More precisely, if $q(x_0)$ is prime for some $x_0$, we can use the Complex Multiplication method \cite{bss}, \cite{mor} to construct an elliptic curve over $\F_{q(x_0)}$ with $n(x_0)$ points and embedding degree $k$.  

\begin{theorem}
\label{t:method}
Fix a positive integer $k$, and let $\Phi_k(x)$ be the $k$th cyclotomic polynomial.  Let $t(x)$ be a polynomial with integer coefficients, let $n(x)$ be an irreducible factor of $\Phi_k(t(x)-1)$, and let $q(x) = n(x) + t(x) - 1$.  Let $f(x) = 4q(x) - t(x)^2$.  Fix a positive square-free integer $D$, and suppose $(x_0,y_0)$ is an integer solution to the equation $Dy^2 = f(x)$ for which 
\begin{enumerate}
	\item $q(x_0)$ is prime, and
	\item $n(x_0)$ is prime.  
\end{enumerate}
If $D$ is sufficiently small, then there is an efficient algorithm to construct an elliptic curve $E$ defined over $\F_{q(x_0)}$ such that $E(\F_{q(x_0)})$ has prime order $n(x_0)$ and $E$ has embedding degree at most $k$.
\end{theorem}

\begin{proof}
By hypothesis, we have a solution $(x_0,y_0)$ to the equation $Dy^2 = f(x)$ for which $q(x_0)$ is prime.  If $D$ is sufficiently small then the construction of an elliptic curve $E$ over $\F_{q(x_0)}$ with $\#E(\F_{q(x_0)}) = n(x_0)$ is standard via the Complex Multiplication method; see \cite{bss} or \cite{mor} for details.  Since $n(x_0)$ is prime, $E(\F_{q(x_0)})$ has prime order, and it remains only to show that $E$ has embedding degree at most $k$.  Barreto, Lynn, and Scott \cite[Lemma 1]{bls1} show that $E$ having embedding degree $k$ is equivalent to $n(x_0)$ dividing $\Phi_k(t(x_0)-1)$ and $n(x_0)$ not dividing $\Phi_i(t(x_i)-1)$ for $i < k$.  Since we have chosen the polynomial $n(x)$ to divide $\Phi_k(t(x)-1)$, $n(x_0)$ is guaranteed to divide $q(x_0)^k-1$, and the embedding degree of $E$ is thus at most $k$. 
\qed
\end{proof}

\begin{remark}
The fact that $n(x)$ does not divide $\Phi_i(t(x)-1)$ as polynomials for $i < k$ does not guarantee that $n(x_0)$ does not divide $\Phi_i(t(x_0)-1)$ as integers for some $i < k$.  However, this latter case will be rare in practice, and thus the embedding degree of a curve constructed via the method of Theorem \ref{t:method} will usually be $k$.
\end{remark}

\begin{remark}
If we wish to construct curves whose orders are not necessarily prime but merely have a large prime factor, we may relax condition (2) of the theorem accordingly, and the same analysis holds.
\end{remark}

In practice, to construct an elliptic curve with embedding degree $k$ one chooses polynomials $t(x)$, $n(x)$, and $q(x)$ satisfying the conditions of Theorem \ref{t:method} and tests various values of $x$ until $n(x)$ and $q(x)$ are prime.  If the distributions of the values of the polynomials $n(x)$ and $q(x)$ are sufficiently random, the Prime Number Theorem tells us that we should have to test roughly $\log n(x_1) \log q(x_1)$ values of $x$ near $x_1$ until we find an $x_0$ that gives a prime value for both polynomials.  Since the distribution of prime values of polynomials is not well understood in general, it will be hard to prove theorems that explicitly construct infinite families of elliptic curves of prime order.  Rather, we will be slightly less ambitious and search for polynomials as in Theorem \ref{t:method} that will give us the desired elliptic curves whenever the polynomials take on prime values.  We incorporate this approach into the following definition.

\begin{definition}
\label{d:family}
Let $t(x)$, $n(x)$, and $q(x)$ be polynomials with integer coefficients.  For a given positive integer $k$ and positive square-free integer $D$, the triple $(t,n,q)$ {\it represents a family of curves with embedding degree $k$} if the following conditions are satisfied:
\begin{enumerate}
\item $n(x) = q(x) + 1 - t(x)$.
\item $n(x)$ and $q(x)$ are irreducible.
\item $n(x)$ divides $\Phi_k(t(x)-1)$, where $\Phi_k$ is the $k$th cyclotomic polynomial.
\item The equation $Dy^2 = 4q(x) - t(x)^2$ has infinitely many integer solutions $(x,y)$.
\end{enumerate}
\end{definition}

Defining a family of curves in this way gives us a simple criterion for constructing elliptic curves with embedding degree $k$.  This criterion is implicit in the Barreto-Naehrig construction of curves with $k = 12$ and $D = 3$ \cite{bn}.

\begin{corollary}
Suppose $(t,n,q)$ represents a family of curves with embedding degree $k$ for some $D$.  Then for each $x_0$ such that $n(x_0)$ and $q(x_0)$ are both prime, there is an elliptic curve $E$ defined over $\F_{q(x_0)}$ such that $\#E(\F_{q(x_0)})$ is prime, and $E$ has embedding degree at most $k$.
\end{corollary}

In practice, for any $t(x)$ we can easily find $n(x)$ and $q(x)$ satisfying conditions (1), (2), and (3) of Definition \ref{d:family}; the difficulty arises in choosing the polynomials so that $Dy^2 = 4q(x) - t(x)^2$ has infinitely many integer solutions.  In general, if $f(x)$ is a square-free polynomial of degree at least $3$, then there will be only a finite number of integer solutions to the equation $Dy^2 = f(x)$ (cf.\ Proposition \ref{p:nofamily}).  Thus we conclude that $(t,n,q)$ can represent a family of curves only if $f(x)$ has some kind of special form.  

We now show that if $f(x)$ is quadratic, then one integral solution to the equation $Dy^2 = f(x)$ will give us infinitely many solutions.  This is the technique that Miyaji, et al.\ \cite{mnt} use to produce curves with embedding degree 3, 4, or 6, and we will use the same technique in Section \ref{s:k=10} to construct curves with embedding degree $10$.  The idea is as follows: we complete the square to write the equation $Dy^2 = f(x)$ as $u^2 - D'v^2 = T$ for some constant $T$, and observe that $(u,v)$ is a solution to this equation if and only if $u + v\sqrt{D'}$ has norm $T$ in the real quadratic field $\Q(\sqrt{D'})$.  By Dirichlet's unit theorem, there is a one-dimensional set of norm-one integral elements of this field; multiplying each of these units by our element of norm $T$ gives an infinite family of elements of norm $T$.  We then show that a certain fraction of these elements can be converted back to solutions of the original equation.  

\begin{theorem}
\label{t:quadratic}
Fix an integer $k > 0$, and choose polynomials $t(x)$, $n(x)$, $q(x) \in \Z[x]$ satisfying conditions (1), (2), and (3) of Definition \ref{d:family}.  Let $f(x) = 4 q(x)- t(x)^2$.  Suppose $f(x) = ax^2+bx+c$, with $a,b,c \in \Z$, $a > 0$, and $b^2 - 4ac \neq 0$.  Let $D$ be a square-free integer such that $aD$ is not a square.  If the equation $Dy^2 = f(x)$ has a solution $(x_0,y_0)$ in the integers, then $(t,n,q)$ represents a family of curves with embedding degree $k$.
\end{theorem}

\begin{proof}
Completing the square in the equation $Dy^2 = f(x)$ and multiplying by $4a$ gives
\begin{equation}
aD(2y)^2 = (2ax+b)^2 - (b^2-4ac).
\end{equation}
If we write $aD = D'r^2$ with $D'$ square-free and make the substitutions $u = 2ax+b$, $v = 2ry$, $T = b^2-4ac$, the equation becomes
\begin{equation}
u^2 - D'v^2 = T.
\end{equation}
Note that since $aD$ is not a square, we have $D' > 1$. 

Under the above substitution, a solution $(x_0,y_0)$ to the original equation $Dy^2 = f(x)$ gives an element $u_0 + v_0 \sqrt{D'}$ of the real quadratic field $\Q(\sqrt{D'})$ with norm $T$.  Furthermore, this solution satisfies the congruence conditions
\begin{equation}
\label{eq:cong}
\begin{split}
u_0 & \equiv  b \pmod{2a} \\
v_0 & \equiv  0 \pmod{2r}.
\end{split}
\end{equation}
We wish to find an infinite set of solutions $(u,v)$ satisfying the same congruence conditions, for we can transform such a solution into an integer solution to the original equation.  To find such solutions we employ Dirichlet's unit theorem \cite[\S 1.7]{neu}, which tells us that the integer solutions to the equation $\alpha^2 - D'\beta^2 = 1$ are in one-to-one correspondence with the real numbers $\alpha + \beta \sqrt{D'} = \pm (\alpha_0 + \beta_0 \sqrt{D'})^n$ for some fixed $(\alpha_0,\beta_0)$ and any integer $n$.  
The real number $\alpha_0 + \beta_0 \sqrt{D'}$ is either a fundamental unit of the real quadratic field $\Q(\sqrt{D'})$ or (if the norm of the fundamental unit is $-1$) the square of a fundamental unit.

Reducing the coefficients of $\alpha_0 + \beta_0 \sqrt{D'}$ modulo $2a$ gives an element $z = \bar \alpha_0 + \bar \beta_0 \bar x$ of the ring 
\begin{equation}
	R = \frac{\Z/{2a\Z}[x]}{(x^2 - D')}.
\end{equation}
Furthermore, since $(\alpha_0 + \beta_0 \sqrt{D'})(\alpha_0 - \beta_0 \sqrt{D'}) = 1$, $z$ is invertible in $R$, i.e.\ $z \in R^*$.  Since $R^*$ is a finite group of size less than $4a^2$, there is an integer $m < 4a^2$ such that $z^m = 1$ in $R^*$.\footnote{In fact, since $z$ is an element of the norm-one subgroup of $R^*$, $m$ is bounded above by $2^sa$, where $s$ is the number of distinct primes dividing $2a$.  A more detailed study of the group $R^*$ appears in an earlier draft of this paper \cite{free}.}  Lifting back up to the full ring $\Z[\sqrt{D'}]$, we see that $(\alpha_0 + \beta_0 \sqrt{D'})^m = \alpha_1 + \beta_1 \sqrt{D'}$ for integers $\alpha_1, \beta_1$ satisfying
\begin{equation}
\label{eq:cong2}
\begin{split}
\alpha_1 & \equiv  1  \pmod{2a}, \\
\beta_1 & \equiv  0 \pmod{2a}.
\end{split}
\end{equation}

Now for any integer $n$ we can compute integers $(u,v)$ such that 
\begin{equation}
u + v\sqrt{D'} = (u_0 + v_0 \sqrt{D'})(\alpha_1 + \beta_1 \sqrt{D'})^n.
\end{equation}
We claim that $(u,v)$ satisfy the congruence conditions (\ref{eq:cong}).
To see this, let $\alpha_n + \beta_n \sqrt{D'} = (\alpha_1 + \beta_1 \sqrt{D'})^n$.  The conditions (\ref{eq:cong2}) imply that $\alpha_n \equiv 1 \pmod{2a}$ and $\beta_n \equiv 0 \pmod{2a}$.  Combining this observation with the formulas
\begin{equation}
\begin{split}
    u & =  \alpha_n u_0 + \beta_n v_0 D' \\
    v & =  \alpha_n v_0 + \beta_n u_0,
\end{split}
\end{equation}
we see that $u \equiv u_0 \equiv b \pmod{2a}$ and $v \equiv v_0 \pmod{2a}$.  Furthermore, $v_0 \equiv 0 \pmod{2r}$ and $2r$ divides $2a$ (since $aD = D'r^2$ and $D$ is square-free), so we conclude that $v \equiv 0 \pmod{2r}$.

The new solution $(u,v)$ thus satisfies the congruence conditions (\ref{eq:cong}).  Any integer $n$ gives such a solution, so by setting $x = (u-b)/2a$ and $y = v/2r$ for each such $(u,v)$, we have generated an infinite number of integer solutions to the equation $Dy^2 = f(x)$.  This is condition (4) of Definition \ref{d:family}; by hypothesis $(t,n,q)$ satisfy conditions (1), (2), and (3), so we conclude that $(t,n,q)$ represents a family of curves with embedding degree $k$.
\qed
\end{proof}

\begin{remark}
More generally, we may find an infinite family of curves in the case where $f(x) = g(x)^2 h(x)$, with $h(x)$ quadratic.  Specifically, if we let $y = y'g(x)$, then given one integral solution $(x,y')$ to the equation $Dy'^2 = h(x)$ we may use the method of Theorem \ref{t:quadratic} to find an infinite number of solutions.  However, we currently know of no examples for which $f(x)$ is of this form.
\end{remark}

Theorem \ref{t:quadratic} tells us that if $f(x)$ is quadratic and square-free, we may get a family of curves of the prescribed embedding degree for {\it each} $D$.  If $f(x)$ is instead a linear function times a square, then we may still get a family of curves, but for only a single $D$.  This is the method that Barreto and Naehrig \cite{bn} use to construct curves with $k = 12$ (see Section \ref{ss:12}).

\begin{proposition}
\label{p:square}
Fix an integer $k > 0$, and let $n(x)$, $t(x)$, and $q(x)$ be polynomials in $\Z[x]$ satisfying conditions (1), (2), and (3) of Definition \ref{d:family}.  Let $f(x) = 4 q(x)- t(x)^2$, and suppose $f(x) = (Ax+D) g(x)^2$ for some positive integer $D$ and some polynomial $g(x)$.  Then $(t,n,q)$ represents a family of curves with embedding degree $k$.
\end{proposition}

\begin{proof}
For any integer $v$, we set $x = ADv^2 + 2Dv$ and let $y = (Av+1)g(x)$.  An easy computation shows that $(x,y)$ is a solution to the equation $Dy^2 = f(x)$, so if $D$ is square-free then condition (4) is satisfied for the integer $D$.  If $D$ is not square-free then we may absorb its square factors into $y$, and condition (4) is satisfied for the largest square-free factor $D'$ of $D$.
\qed
\end{proof}

We conclude this section with a partial converse to Theorem \ref{t:quadratic}; namely, if the degree of $f(x)$ is at least $3$, then we are unlikely to find an infinite family of curves.

\begin{proposition}
\label{p:nofamily}
Let $(t,n,q)$ be polynomials with integer coefficients satisfying conditions (1), (2), and (3) of Definition \ref{d:family}, and let $f(x) = 4q(x) - t(x)^2$.  Suppose $f(x)$ is square-free and $\deg f(x) \geq 3$.  Then $(t,n,q)$ does not represent a family of elliptic curves with embedding degree $k$.
\end{proposition}

\begin{proof}
Since $f(x)$ is square-free (i.e.\ has no double roots) and has degree at least $3$, the equation $Dy^2 = f(x)$ defines a smooth affine plane curve of genus $g \geq 1$.  By Siegel's Theorem (cf.\ \cite[Theorem IX.4.3]{sil} and \cite[\S I.2]{cs}) such curves have a finite number of integral points, so condition (4) is not satisfied.
\qed
\end{proof}

\section{Elliptic Curves with Embedding Degree $10$.}
\label{s:k=10}

In this section, we use the method of Section \ref{s:method}, and Theorem \ref{t:quadratic} in particular, to construct elliptic curves of prime order with embedding degree $10$.  Our key observation is that since the hypotheses of Theorem \ref{t:quadratic} require $f(x) = 4n(x) - (t(x)-2)^2$ to be quadratic, we should choose $n(x)$ and $t(x)$ in such a way that the high-degree terms of $t(x)^2$ cancel out those of $4n(x)$; in particular, the degree of $t(x)$ must be half the degree of $n(x)$.  We have discovered that for $k = 10$ there is a choice of $n(x)$ and $t(x)$ such that this is possible.  The resulting construction of elliptic curves with embedding degree $10$ solves an open problem posed by Boneh, Lynn, and Shacham \cite[\S 4.5]{bls}.

We begin by recalling that to construct a curve with embedding degree $k$, we must choose the number of points $n(x)$ and the trace $t(x)$ such that $n(x)$ is an irreducible factor of $\Phi_k(t(x) - 1)$, where $\Phi_k$ is the $k$th cyclotomic polynomial.  If $k = 10$ and $t(x)$ is linear then $\Phi_k(t(x)-1)$ is an irreducible quartic polynomial, so there is no hope of $f(x) = 4n(x) - (t(x)-2)^2$ being quadratic.  If $k = 10$ and $t(x)$ is quadratic,  Galbraith, McKee, and Valen\c ca \cite{gmv} show that in this case $\Phi_k(t(x)-1)$ either is irreducible of degree $8$ or factors into two irreducible quartic polynomials.  They then show that there is an infinite set of $t(x)$ such that the latter occurs, and that these $t(x)$ are parameterized by the rational points of a certain elliptic curve.  By experimenting with some of the examples given by Galbraith, et al., we discovered that $t(x) = 10x^2 + 5x + 3$ leads to a quadratic $f(x)$.

\begin{theorem}
\label{t:embed10}
Fix a positive square-free integer $D$ relatively prime to $15$.  Define $t(x)$, $n(x)$, and $q(x)$ by
\begin{eqnarray*}
    t(x) & = & 10x^2 + 5x + 3 \\
    n(x) & = & 25x^4+25x^3+15x^2+5x+1 \\
    q(x) & = & 25x^4 + 25x^3 + 25x^2 + 10x + 3.
\end{eqnarray*}
If the equation $u^2 - 15Dv^2 = -20$ has a solution with $u \equiv 5 \pmod{15}$, then $(t,n,q)$ represents a family of curves with embedding degree $10$.
\end{theorem}

\begin{proof}
It is easy to verify that conditions (1)-(3) of Definition \ref{d:family} hold.  Condition (4) requires an infinite number of integer solutions to $Dy^2 = f(x)$, where $f(x) = 4q(x) - t(x)^2$.  The key observation is that for this choice of $t$ and $n$,
\begin{equation}
f(x) = 4q(x) - t(x)^2 = 15x^2 + 10x + 3.
\end{equation}
Multiplying by $15$ and completing the square transforms the equation we wish to solve into
\begin{equation}
D'y^2 = (15x + 5)^2 + 20,
\end{equation}
where $D' = 15D$.
Integer solutions to this equation correspond to integer solutions to $u^2 - D'v^2 = -20$ with $u \equiv 5 \pmod{15}$.  By Theorem \ref{t:quadratic}, if one such solution exists then an infinite number exist, so $(t,n,q)$ represents a family of curves with embedding degree $10$.
\qed
\end{proof}

To use the above result to construct curves with embedding degree $10$, we choose a $D$ and search for solutions to the equation $u^2 - 15Dv^2 = -20$ that give prime values for $q$ and $n$.  The following lemma, proposed by Mike Scott, speeds up this process by restricting the values of $D$ that we can use.

\begin{lemma}
\label{l:cong}
Let $q(x)$ be as in Theorem \ref{t:embed10}.  If $(x,y)$ is an integer solution to $Dy^2 = 15x^2 + 10x + 3$ such that $q(x)$ is prime, then $D \equiv 43$ or $67 \pmod{120}$.
\end{lemma}

\begin{proof}
If $x \equiv 0$ or $2 \pmod 3$ then $q(x)$ is divisible by $3$, while if $x$ is odd then $q(x)$ is even.  Thus if $q(x)$ is prime, then $x \equiv 4 \pmod{6}$.

To deduce the stated congruence for $D$, we consider the equation $Dy^2 = 15x^2 + 10x + 3$ modulo $3$, $5$, and $8$.  To begin, we have $Dy^2 \equiv x \equiv 1 \pmod{3}$, so $D \equiv 1 \pmod{3}$.  Next, we have $Dy^2 \equiv 3 \pmod{5}$, so $y^2 \equiv 1$ or $4 \pmod{5}$ and $D \equiv 2$ or $3 \pmod{5}$.  Finally, since $x$ is even we see that $Dy^2 = 3 \pmod{8}$, and thus $y^2 \equiv 1 \pmod{8}$ and $D \equiv 3 \pmod{8}$.  Combining these results via the Chinese remainder theorem, we conclude that $D \equiv 43$ or $67 \pmod{120}$.
\qed
\end{proof}

After reading an earlier draft of this paper \cite{free}, Mike Scott used Theorem \ref{t:embed10} and Lemma \ref{l:cong} to find examples of elliptic curves with embedding degree $10$ via the following algorithm.

\begin{enumerate}
	\item Choose a $D$ such that $15D$ is square-free and $D \equiv 43$ or $67 \pmod{120}$.
	\item Find solutions $(u,v)$ to the equation $u^2 - 15Dv^2 = -20$.
	\item For each solution $(u,v)$:
	\begin{enumerate}
		\item If $u$ is too large (e.g.\ $\geq 128$ bits), go to the next solution.
		\item If $u \equiv \pm 5 \pmod{15}$, then
		\begin{enumerate}
			\item Let $x = (-5 \pm u)/15$.
 			\item If $q(x)$ and $n(x)$ are prime, output $(D,x)$
		\end{enumerate}
		\item Multiply $u + v \sqrt{15D}$ by a norm-one element of $\Q(\sqrt{15D})$ to get a new $u$, and return to step (a).
	\end{enumerate}
	\item Increase $D$ and return to step (1).
\end{enumerate}

For each $(D,x)$ output by the algorithm, Scott used the Complex Multiplication method (cf.\ \cite{bss}, \cite{mor}) to construct an elliptic curve over $\F_{q(x)}$ whose number of points is $n(x)$.  By Theorem \ref{t:method} this curve has embedding degree at most $10$, and in practice we find that the embedding degree is exactly $10$.  Below are two examples of elliptic curves constructed in this manner.  

\begin{example} 
\label{ex1}
(A 149-bit curve.) Choosing $D = 1666603$ and running the above algorithm produces the following example.  Let $q, n, A, B$ be as follows:
\begin{eqnarray*}
q & = & 503189899097385532598615948567975432740967203 \\
n & = & 503189899097385532598571084778608176410973351 \\
A & = & -3 \\
B & = & 78778770898368212452154728282767760988008151.
\end{eqnarray*}
Then $q$ and $n$ are $149$-bit prime numbers such that the curve $y^2 = x^3 + Ax + B$ defined over $\F_q$ has $n$ points.  Since $n \mid q^{10} - 1$ and $n \nmid q^i - 1$ for $i < 10$, this curve has embedding degree $10$.  
\end{example}

\begin{example} 
\label{ex2} 
(A 196-bit curve.) Choosing $D = 579003643$ and running the above algorithm produces the following example. Let $q, n, A, B$ be as follows:
\begin{eqnarray*}
q & = & 61099963271083128746073769567944870354270161646150914794603 \\
n & = & 61099963271083128746073769567450502219087145916434839626301 \\
A & = & -3 \\
B & = & 1112775869471458154129950648198203893613615552476491488167.
\end{eqnarray*}
Then $q$ and $n$ are $196$-bit prime numbers such that the curve $y^2 = x^3 + Ax + B$ defined over $\F_q$ has $n$ points.  Since $n \mid q^{10} - 1$ and $n \nmid q^i - 1$ for $i < 10$, this curve has embedding degree $10$.  
\end{example}

\section{Elliptic Curve Families with Small Embedding Degree}
\label{s:small}

In this section we show how the existing constructions of ordinary elliptic curves of prime order with embedding degree 3, 4, or 6 \cite{mnt} or embedding degree 12 \cite{bn} can be explained via the  framework of Section \ref{s:method}.  The former uses Theorem \ref{t:quadratic}, while the latter employs Proposition \ref{p:square}.

\subsection{MNT Elliptic Curves}
\label{ss:mnt}

Miyaji, Nakabayashi, and Takano \cite{mnt} have classified all ordinary elliptic curves of prime order with embedding degree $3$, $4$, and $6$.  Their theorem is as follows:

\begin{theorem}[\cite{mnt}]  Let $E$ be an ordinary elliptic curve over $\F_q$ such that $\#E(\F_q) = n = q+1-t$ is prime and $E$ has embedding degree $k = 3,4$, or $6$.  Then there exists an integer $x$ such that $t$, $n$, and $q$ are of the form specified in the following table:
$$\begin{array}{|c|c|c|c|}
\hline
k & t & n & q \\
\hline
3 & -1 \pm 6x & 12x^2 \mp 6x + 1 & 12x^2 - 1 \\
4 & -x \mbox{ or } x+1 & x^2 + 2x + 2 \mbox{ or } x^2 + 1 &  x^2+x+1 \\
6 & 1 \pm 2x & 4x^2 \mp 2x + 1 & 4x^2 + 1 \\
\hline
\end{array} $$
\end{theorem}

This theorem fits into the framework of Section \ref{s:method} as follows.  To find an infinite family of curves via Theorem \ref{t:quadratic}, we require $f(x)$ to be quadratic.  Since $\deg \Phi_k(x) = 2$ for $k = 3$, $4$, or $6$, if we let $t(x)$ be any linear polynomial and $n(x)$ be the (irreducible) quadratic $\Phi_k(t(x)-1)$ (with any constant factor divided out), then $f(x) = 4n(x) - (t(x)-2)^2$ is quadratic.  If $q(x) = n(x) + t(x) - 1$ is also irreducible and the equation $Dy^2 = f(x)$ has one solution, then $(t,n,q)$ satisfy the hypotheses of Theorem \ref{t:quadratic} and thus represent a family of curves with embedding degree $k$.  Miyaji, et al.\ arrive at their stronger result by using the fact that $\#E(\F_q)$ is prime to show that any values of $t$, $n$, and $q$ that give rise to such a curve must be of the specified form.

\subsection{Elliptic Curves with Embedding Degree 12}
\label{ss:12}

Finally, we note that the Barreto-Naehrig construction \cite{bn} of curves with embedding degree $12$ falls under the case of Proposition \ref{p:square}.  Specifically, if $t(x) = 6x^2 + 1$, then $\Phi_{12}(t(x) - 1) = n(x)n(-x)$, where $n(x) = 36x^4 + 36x^3 + 18x^2 + 6x + 1$, and  
\begin{equation}
f(x) = 4n(x) - (t(x) - 2)^2 = 3(6x^2 + 4x + 1)^2.
\end{equation}
Since $q(x) = 36x^4 + 36x^3 + 12x^2 + 6x + 1$ is also irreducible, Proposition \ref{p:square} tells us that if we set $D = 3$, then $(t,n,q)$ represents a family of curves with embedding degree $12$.

\section{Higher Embedding Degrees}
\label{s:higher}

To construct families of elliptic curves with prescribed embedding degree, the method of Section \ref{s:method} requires us to find an infinite number of integer solutions to an equation of the form $Dy^2 = f(x)$.  In this section, we give evidence that in general the degree of $f(x)$ is large, and thus by Proposition \ref{p:nofamily} we are unlikely to find an infinite family of curves.  We begin with a lemma that restricts the possible degrees of the polynomial $n(x)$; the lemma generalizes a result of Galbraith, et al.\ \cite[Lemma 1]{gmv}.

\begin{lemma}
\label{l:degn}
Fix $k$, let $t(x)$ be a polynomial, and let $n(x)$ be an irreducible factor of $\Phi_k(t(x)-1)$.  Then the degree of $n$ is a multiple of $\varphi(k)$, where $\varphi$ is the Euler phi function.
\end{lemma}

\begin{proof}
Suppose $t(x)$ has degree $d$, so $\deg \Phi_k(t(x)-1) = d \varphi(k)$.  Let $\theta$ be a root of $n(x)$, and let $\omega = t(\theta)-1$.  Then $\Phi_k(\omega) = 0$, so $\omega$ is a primitive $k$th root of unity.  We thus have the inclusion of fields $\Q(\theta) \supset \Q(\omega) \supset \Q$.  Since $[\Q(\theta):\Q] = \deg n(x)$ and $[\Q(\omega):\Q] = \varphi(k)$, we conclude that $\varphi(k)$ divides $\deg n(x)$.
\qed
\end{proof}

The key observation that allowed us to construct families of elliptic curves with embedding degree $10$ was that if $f(x)$ is quadratic and $n(x)$ has degree greater than $2$, then the polynomial $t(x)$ must be chosen so that the high degree terms of $t(x)^2$ cancel out those of $4n(x)$.  The following proposition shows that this is in fact the only way to construct such families.

\begin{proposition}
\label{p:higher}
Suppose $(t,n,q)$ represents a family of curves with embedding degree $k$, and suppose further that $f(x) = 4n(x) - (t(x)-2)^2$ is square-free.  If $\varphi(k) \geq 4$, then 
\begin{equation}
\deg t(x) =  \frac{1}{2} \deg n(x) = \frac{1}{2} \deg q(x).
\end{equation}
Furthermore, if $a$ is the leading coefficient of $t(x)$, then $a^2/4$ is the leading coefficient of $n(x)$ and $q(x)$.
\end{proposition}

\begin{proof}
Since $\varphi(k) \geq 4$, by Lemma \ref{l:degn} $\deg n(x) \geq 4$, and since $f(x)$ is square-free, by Proposition \ref{p:nofamily} $\deg f(x) \leq 2$.  Since $f(x) = 4n(x) - (t(x)-2)^2$, we conclude that $\deg t(x) = \frac{1}{2} \deg n(x)$, and since $n(x) = q(x) + 1 - t(x)$, we see that $\deg n(x) = \deg q(x)$.  The observation about the leading coefficients follows immediately.
\qed
\end{proof}

As an immediate corollary, we see that if $k > 6$ (so $\varphi(k) \geq 4$) then choosing a linear $t(x)$ will not in general give us an infinite family of curves, whereas if $k > 12$ (so $\varphi(k) \geq 6$) then choosing a quadratic $t(x)$ will not in general give us an infinite family of curves.  

Proposition \ref{p:higher} tells us that for embedding degrees $k$ with $\varphi(k) \geq 4$, to find an infinite family of curves we will have to choose $t(x_0)$ of degree at least $2$ such that $\phi_k(t(x)-1)$ is not irreducible.  Galbraith, McKee, and Valen\c ca \cite{gmv} observe that this is hard even for quadratic $t(x)$, and as the degree increases the problem will only become more difficult.  An alternative would be to choose $t$ and $n$ such that $f(x)$ has a square factor; this appears to be just as difficult, but has not been studied in depth.

\section {Conclusion}
\label{s:conclusion}

We have seen in Section \ref{s:method} that the current methods for constructing families of elliptic curves of prime order with prescribed embedding degree can all be subsumed under a general framework.  In Section \ref{s:k=10} we showed how this framework can be used to construct curves with embedding degree $10$ and we gave examples of such curves, which have not previously appeared in the literature.  In Section \ref{s:small} we showed how this framework incorporates the existing constructions for embedding degrees $3,4,6$, and $12$.  

In Section \ref{s:higher} we showed that our method can only produce an infinite family of curves if a certain polynomial $f(x)$ either is quadratic or has a square factor.  These two conditions have been achieved for $k = 10$ and $k = 12$, respectively, but these two examples appear to be special cases, and in general we have not found a way to achieve either of these two conditions.  The success of our method in producing curves with embedding degree greater than $12$ depends on our ability to control the behavior of $f(x)$, which leads to the following important open problem.

\begin{problem}
\label{q:polys}
Given an integer $k$ such that $\varphi(k) \geq 4$, find polynomials $t(x)$ and $n(x)$ such that
\begin{enumerate}
\item $n(x)$ is an irreducible factor of $\Phi_k(t(x)-1)$, where $\Phi_k$ is the $k$th cyclotomic polynomial, and
\item $f(x) = 4n(x) - (t(x)-2)^2$ is either quadratic or of the form $g(x)^2 h(x)$, with $\deg h(x) \leq 2$.
\end{enumerate}
\end{problem}


\begin{thebibliography}{MMM}

\bibitem{bls1} P.S.L.M. Barreto, B. Lynn, M. Scott, ``Constructing elliptic curves with prescribed embedding degrees,'' in {\it SCN 2002}, ed. S. Cimato, C. Galdi, G. Persiano, Springer LNCS {\bf 2576} (2003) 257-267.

\bibitem{bn} P.S.L.M. Barreto, M. Naehrig, ``Pairing-friendly elliptic curves of prime order," Cryptology ePrint Archive, report 2005/133 (2005), available at {\tt http://eprint.iacr.org/2005/133}.

\bibitem{bss} I. Blake, G. Seroussi, N. Smart, {\it Elliptic Curves in Cryptography}, LMS Lecture Note Series {\bf 265}, Cambridge University Press, 1999.

\bibitem{bss2} I. Blake, G. Seroussi, N. Smart, eds., {\it Advances in Elliptic Curve Cryptography}, LMS Lecture Note Series {\bf 317}, Cambridge Unviersity Press, 2005.

\bibitem{bf} D. Boneh, M. Franklin, ``Identity based encryption from the Weil pairing,'' in {\it CRYPTO '01}, ed.\ J. Kilian, Springer LNCS {\bf 2139} (2001), 213-229.

\bibitem{bls} D. Boneh, B. Lynn, H. Shacham, ``Short signatures from the Weil pairing,'' in {\it ASIACRYPT '01}, ed. C. Boyd, Springer LNCS {\bf 2248} (2001), 514-532.

\bibitem{bw} F. Brezing, A. Weng, ``Elliptic curves suitable for pairing based cryptography,'' to appear in {\it Designs, Codes, and Cryptography}; preprint available at {\tt http://eprint.iacr.org/2003/143}.

\bibitem{cp} C. Cocks, R.G.E. Pinch, ``Identity-based cryptosystems based on the Weil pairing," unpublished manuscript, 2001.

\bibitem{cs} G. Cornell, J. Silverman, eds., {\it Arithmetic Geometry,} Springer, New York 1986.

\bibitem{cdc} S. Cui, P. Duan, C.W. Chan, ``A new method of building more non-supersingular elliptic curves,'' in {\it ISH 2005}, ed. O. Gervasi et al., Springer LNCS {\bf 3481} (2005), 657-664.

\bibitem{free} D. Freeman, ``Constructing families of pairing-friendly elliptic curves,'' Hew-lett-Packard Laboratories technical report HPL-2005-155 (2005), available at {\tt http://www.hpl.hp.com/techreports/2005/HPL-2005-155.html}.

\bibitem{gmv} S. Galbraith, J. McKee, P. Valen\c ca, ``Ordinary abelian varieties having small embedding degree," Cryptology ePrint Archive, report 2004/365 (2004), available at {\tt http://eprint.iacr.org/2004/365}.

\bibitem{joux} A. Joux, ``A one round protocol for tripartite Diffie-Hellman,'' in {\it ANTS-IV}, ed. W. Bosma, Springer LNCS {\bf 1838} (2000), 385-394.

\bibitem{mov} A. Menezes, T. Okamoto, S. Vanstone, ``Reducing elliptic curve logarithms to logarithms in a finite field," {\it IEEE Transactions on Information Theory} {\bf 39} (1993), 1639-1646.

\bibitem{mnt} A. Miyaji, M. Nakabayashi, S. Takano, ``New explicit conditions of elliptic curve traces for FR-reduction," {\it IEICE Transactions on Fundamentals} {\bf E84-A(5)} (2001), 1234-1243.

\bibitem{mor} F. Morain, ``Building cyclic elliptic curves modulo large primes," in {\it EUROCRYPT '91}, ed.\ D. W.  Davies, Springer LNCS {\bf 547} (1991) 328-336.

\bibitem{neu} J. Neukirch, {\it Algebraic Number Theory}, Springer, Berlin 1999.

\bibitem{sil} J. Silverman, {\it The Arithmetic of Elliptic Curves}, Springer GTM {\bf 106}, 1986.

\end{thebibliography}
\end{document}